\documentclass[12pt]{article}
\usepackage{epsfig}
\usepackage{epstopdf}
\usepackage{graphicx}
\usepackage{amssymb, amsmath, amsthm}

\newtheorem{dfn}{Definition}
\newtheorem{defn}[dfn]{Definition}

\newtheorem{rem}[dfn]{Remark}

\newtheorem{thm}{Theorem}
\newtheorem{lem}{Lemma}

\newtheorem{prop}[thm]{Proposition}

\newtheorem{cor}[thm]{Corollary}

\def\proof{\par\medskip\noindent{\it Proof: }}

\def\R{{\mathbb R}}

\def\M{{\mathcal M}}

\def\Si{\Sigma}

\def\L{{\cal L}}

\def\>{\rangle}
\def\<{\langle}

\def\D{\partial}
\def\3{\ss}
\def\8{\infty}

\usepackage{lpic}

\begin{document}


\title{The Kakimizu complex is simply connected}
\author{Jennifer Schultens}

\maketitle

\begin{abstract}
  In 1992, Osamu Kakimizu defined a complex that has become known as
  the Kakimizu complex of a knot.  Vertices correspond to isotopy
  classes of minimal genus Seifert surfaces of the knot.  Higher
  dimensional simplices correspond to collections of such classes of
  Seifert surfaces that admit disjoint representatives.  We show that
  this complex is simply connected.

\end{abstract}

\vspace{2 mm}

One of the fundamental objects considered in the topological study of
knots is the Seifert surface.  Interestingly, a knot can have many,
and in some cases infinitely many, non-isotopic Seifert surfaces of minimal genus.  
For instance, a Seifert surface of a connect sum of knots, $K = K_1 \#
K_2$, can be spun around a swallow-follow torus of $K$ to yield
infinitely many Seifert surfaces that are typically non-isotopic.  This phenomenon of
non-uniqueness for Seifert surfaces was described early on in
\cite{Al} and \cite{E}.  The Kakimizu complex aims to capture
structural information of the set of isotopy classes of Seifert
surfaces of a given knot.  It is one of several complexes defined by
considering isotopy classes of certain submanifolds and disjointness
properties of representatives of such isotopy classes.

In 1992, Osamu Kakimizu defined a complex that has become known as
the Kakimizu complex of a knot.  See \cite{K}.  The results in
\cite{Al} and \cite{E} establish that the Kakimizu complex is
nontrivial.  A result of M. Scharlemann and A. Thompson, see
\cite{ST}, establishes that the Kakimizu complex is connected.  It
is conjectured that the Kakimizu complex is contractible.

Distance in the Kakimizu complex can be defined in terms of the
number of edges in an edge path between the two vertices, but
Kakimizu showed that this is equivalent to a more sophisticated
formulation in terms of the universal abelian cover of the knot.
Whereas the first definition is the standard definition for
distance in a complex, the second provides a more effective means
of computing the distance in the Kakimizu complex.

Recent years have seen progress in understanding key facts about
the Kakimizu complex: W. Jaco and E. Sedgwick showed that the
Kakimizu complex of the knot $K$ is finite if $K$ is atoroidal and
has genus at least 2.  See \cite{Oe}.  Moreover, R. Wilson showed
that, in fact, it suffices to assume that $K$ is atoroidal.  See
\cite{W}.  Results pertaining to specific classes of knots can be
found in \cite{HS}, \cite{Pel}, \cite{S}, \cite{SSh} and \cite{Ts}.  
In \cite{SSh}, M. Sakuma and K. Shackleton establish concrete 
diameter bounds and provide an overview of the current understanding 
of the Kakimizu complex.  In particular, they prove that the Kakimizu 
complex is simply connected for knots of genus 1.  A more general
understanding of the shape of the Kakimizu complex is highly
desirable.  Many questions remain unanswered.  Though we establish
simple connectivity here, see Theorem \ref{kaksc}, 
the conjectured contractibility has yet
to be proved.

In Section 1 we provide the formal definition of the Kakimizu complex of a knot
and the notion of distance in the Kakimizu complex.  Section 2 introduces
the concept of a relative least area surface and states two required results.
(These results are proved in the appendix.)  The heart of the paper lies
in Section 3, where we prove two key lemmas that yield information about
weighted paths in the Kakimizu complex.  Section 4 contains the observation
that the Kakimizu complex is a flag complex.  In Section 5 we prove the main
theorem, Theorem \ref{kaksc}, stating that the Kakimizu complex is
simply connected.  In Section 6 we prove that the Kakimizu complex is
contractible in the special case when it is 2-dimensional.  
We finish with a few remarks in Section 7.  The Appendix,
by Misha Kapovich, contains the proofs of the theorems
about relative least area surfaces required in this context and is of independent
interest.  

I wish to thank Jesse Johnson, Makoto Sakuma and Ken Shackleton for
pointing out a mistake in an earlier argument pertaining to Theorem
\ref{kaksc}.  I also wish to thank Misha Kapovich for helpful
conversations and for providing the appendix to this paper.  
This work was supported, in part, by a grant from the
NSF.  It was begun at the Max Planck Institute for Mathematics in
the Sciences located in Leipzig, Germany and completed at the
Max Planck Institute for Mathematics located in Bonn, Germany.  
I wish to thank the institutes for their hospitality.

\section{Preliminaries} \label{prelim}

For basic definitions
concerning knots, see \cite{A}, \cite{L} and \cite{R}.   
For basic definitions concerning complexes, see \cite{BH}.  
For a knot $K$ in ${\mathbb S}^3$ we will denote an open regular
neighborhood of $K$ by $\eta(K)$ and the exterior of $K$,
${\mathbb S}^3 - \eta(K)$, by $E(K)$.  
Recall that a {\em Seifert surface} for $K$ is a connected
surface with connected boundary representing a generator of
$H_2(E(K), \partial E(K))$.

\begin{defn}
  The minimal genus Seifert surfaces of a knot $K$ representing a
  fixed generator of $H_2(E(K), \partial E(K))$ form a simplicial complex as
  follows: 1) Vertices correspond to isotopy classes of minimal genus
  Seifert surfaces; 2) Edges correspond to pairs of vertices admitting
  disjoint representatives and, more generally, $n$-dimensional simplices
  correspond to $(n+1)$-tuples $(v_1, \dots, v_{n+1})$ of vertices admitting
  representatives $(S_1, \dots, S_{n+1})$ such that $S_i \cap S_j = \emptyset$
  for all $i < j$.

This complex is called the {\em Kakimizu complex} of $K$. 
\end{defn}

In our discussion here, paths and loops in a simplicial complex
will traverse only vertices and edges (not higher dimensional
faces).

\begin{defn}
The distance between two vertices $v, v^*$ in the Kakimizu complex
of a knot $K$,
denoted by $d_K(v, v^*)$, is the minimal number of edges in a path
connecting the two vertices.  The length of a loop is the number of
edges in the loop.
\end{defn}

One of the fundamental results concerning the Kakimizu complex is due
to M. Scharlemann and A. Thompson.  See \cite{ST}.  It refers to the intersection
number of surfaces.  Recall that the intersection number of a pair of
surfaces $(S, S^*)$, $i(S, S^*)$, is defined to be the least number of
components of intersection of pairs of surfaces isotopic to $(S, S^*)$
that have transverse intersection. In the language here, the theorem
can be formulated as follows:

\begin{thm} \label{SchTh} (Scharlemann-Thompson) The Kakimizu complex 
of a knot $K$ in ${\mathbb S}^3$ is connected.  Moreover, given two
Seifert surfaces $S, S^*$, the distance of the corresponding vertices
in the Kakimizu complex is bounded above by $i(S, S^*) + 1$.
\end{thm}

\section{A few facts about least area surfaces}

Our arguments will rely extensively on the use of (analytic) least area surfaces.  
In addition, we will be interested in the behavior of our surfaces near the boundary
of our knot complements. 

\begin{defn}
Let $M$ be a compact irreducible smooth manifold with boundary.  A {\em relative} Riemannian metric on $M$ is
a Riemannian metric such that $\partial M$ is strictly convex.  Suppose that $\partial M$ consists of tori and
let ${\cal J}$ be a smooth foliation of $\partial M$ by closed curves.  A properly embedded
surface $F$ in $M$ is {\em relative least area} if the following hold: 1) $\partial F \subset {\cal J}$; 
2) The surface $F$ minimizes area over all surfaces in its proper isotopy class subject to the
constraint $\partial F \subset {\cal J}$.
\end{defn}

Let $M$ be a compact irreducible smooth manifold with $\partial M$ consisting of tori and
with a relative Riemannian metric.  In what follows we will always assume that a foliation
of $\partial M$ is fixed.  In the case of a knot complement, we will assume that ${\cal J}$ consists
of preferred longitudes.  When we consider relative least area surfaces, they will be considered
with respect to this fixed foliation.  Let $F$ be a properly embedded surface in $M$.  
We denote the proper isotopy class of $F$ by $[F]$.  Furthermore, we denote the area of $F$ 
by $A(F)$ and the area of a relative least area representative of $[F]$ by $A([F])$.  

In what follows we will specify a path by the vertices it traverses,
{\it e.g.}, $v^1, \dots, v^n$.

\begin{defn} \label{complexity}
The {\em complexity} of a path $v^1, \dots, v^n$ in the Kakimizu complex is
the ordered pair $(n, a)$ where \[a = A([S^{1}]) + \dots + A([S^n])\]
and $S^{1}, \dots, S^n$ are representatives of
$v^1, \dots, v^n$ respectively. We give the set of complexities the lexicographic order. 

In the special case where $n = 1$, we denote the complexity of the vertex $v^1$ by
$c(v^1)$.
\end{defn}

The following theorems are proved in the appendix.

\begin{thm} \label{laexist}
Let $E$ be a compact irreducible 3-manifold endowed
with a relative Riemannian metric and let $F$ be a properly embedded compact
incompressible surface in $E$. Then there exists a relative least area
representative in the proper isotopy class of $F$.
\end{thm}

This is Corollary \ref{exists} from the appendix.

\begin{thm} \label{min}
If there is a homotopically nontrivial loop of length $n > 0$ in the Kakimizu complex 
of a knot $K$, then there is a homotopically
nontrivial loop of smallest complexity of length $n$.   
\end{thm}

This is Corollary \ref{mini} in the appendix.

\begin{thm} \label{ladisj}
Let $E$ be a compact irreducible 3-manifold endowed
with a relative Riemannian metric and let $F_1, F_2$ be properly embedded
relative least
area incompressible surfaces in $E$ with disjoint representatives
in their proper isotopy classes. Then either $F_1, F_2$ are disjoint
or they coincide.
\end{thm}

This is Theorem \ref{disjoint} in the appendix.
(It is a variant of \cite[Theorem 6.2]{FHS}.)

\begin{figure}[htbp] \centering
 \includegraphics[width=4in]{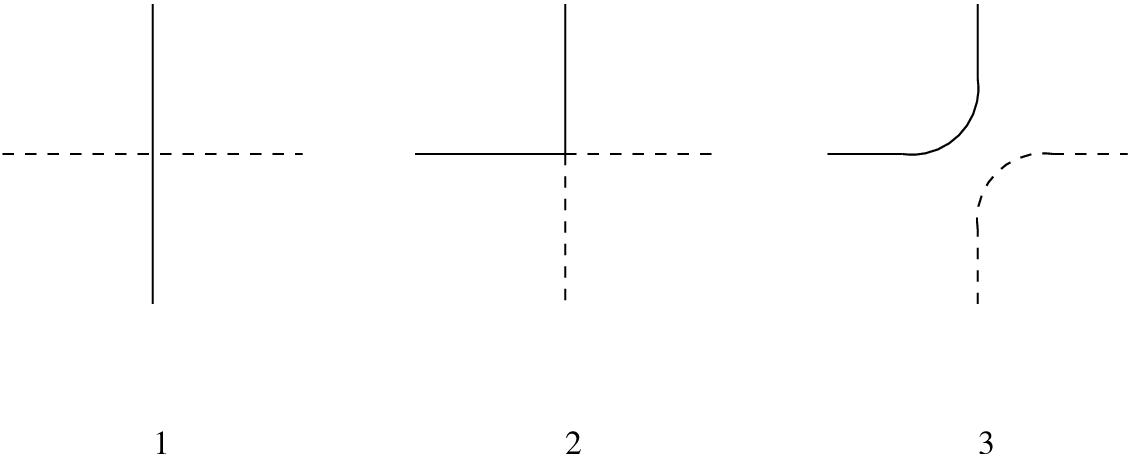} 
\caption{\sl (1) Cross-section of surfaces before exchange. (2) Cross-section of surfaces after exchange. 
(3) Cross-section of surfaces after roundoff.}
\label{xr}
\end{figure}

Two of the standard tools used in conjunction with least area 
surfaces are ``exchange-roundoff"
and the Meeks-Yau trick.  The term ``exchange-roundoff" refers to the fact that
cut-and-paste along a pair of transverse 
least area surfaces yields a pair of lower area least area surfaces.  
(The area remains the same after cut-and-paste, but decreases after roundoff.)  See
Figure \ref{xr}.  

A pair of least area surfaces need not be transverse.  For instance, there can be saddle intersections. 
For specific examples, see \cite{FHS}. Additionally, a pair of relative least area surfaces can share components 
of their boundary.  The Meeks-Yau trick allows us to skirt this issue.  
For an illustration, see the discussion in Case 2 of 
the proof of Lemma \ref{disj} below.  

\section{Key Lemmas}

The following lemmas are crucial in the proof of the main theorem 
(Theorem \ref{kaksc}).  They allow us to refine the
construction of Scharlemann and Thompson \cite{ST}.  We denote
the {\em symmetric difference} of two sets, $X, Y$, by $X \Delta Y$.  Recall
that $X \Delta Y = (X \cup Y) \backslash (X \cap Y)$.  

\begin{lem} \label{disj1}
Let $K$ be a knot in ${\mathbb S}^3$ and
suppose that $S, S^{+1}, S^{-1}$ are minimal genus Seifert surfaces 
for $K$ such that the following hold:

\begin{enumerate}
\item $S$ is disjoint from $S^{+1} \cup S^{-1}$;
\item $S^{-1} \cap S^{+1} \neq \emptyset$;
\item The intersection between  $S^{-1}$ and $S^{+1}$ is transverse and 
$\partial S^{-1} \cap \partial S^{+1}=\emptyset$;
\item There are no disk components in $S^{-1} \Delta S^{+1}$.
\end{enumerate}

\noindent
Then there are two minimal genus Seifert surfaces $S^{up}, S^{down}$ for $K$ such
that (a) \[S \cap S^{up} = \emptyset, S \cap S^{down} = \emptyset, 
(S^{+1} \cup S^{-1}) \cap S^{up} = \emptyset, \] \[
(S^{+1} \cup S^{-1}) \cap S^{down} = \emptyset, S^{up} \cap S^{down} = \emptyset.\]
 
 \noindent
Moreover, (b) if there is a surface $F$ disjoint from $S^{+1} \cup S^{-1}$, 
then it is also disjoint from $S^{down}, S^{up}$.
 
Suppose furthermore that $E(K)$ is endowed with a relative Riemannian metric
and that $S^{+1}, S^{-1}$ are relative least area surfaces.  Then (c)
\[A(S^{-1}) + A(S^{+1}) >  A([S^{down}]) + A([S^{up}]).\] 
\end{lem}

\begin{figure}[htbp] 
   \centering
\includegraphics[width=3in]{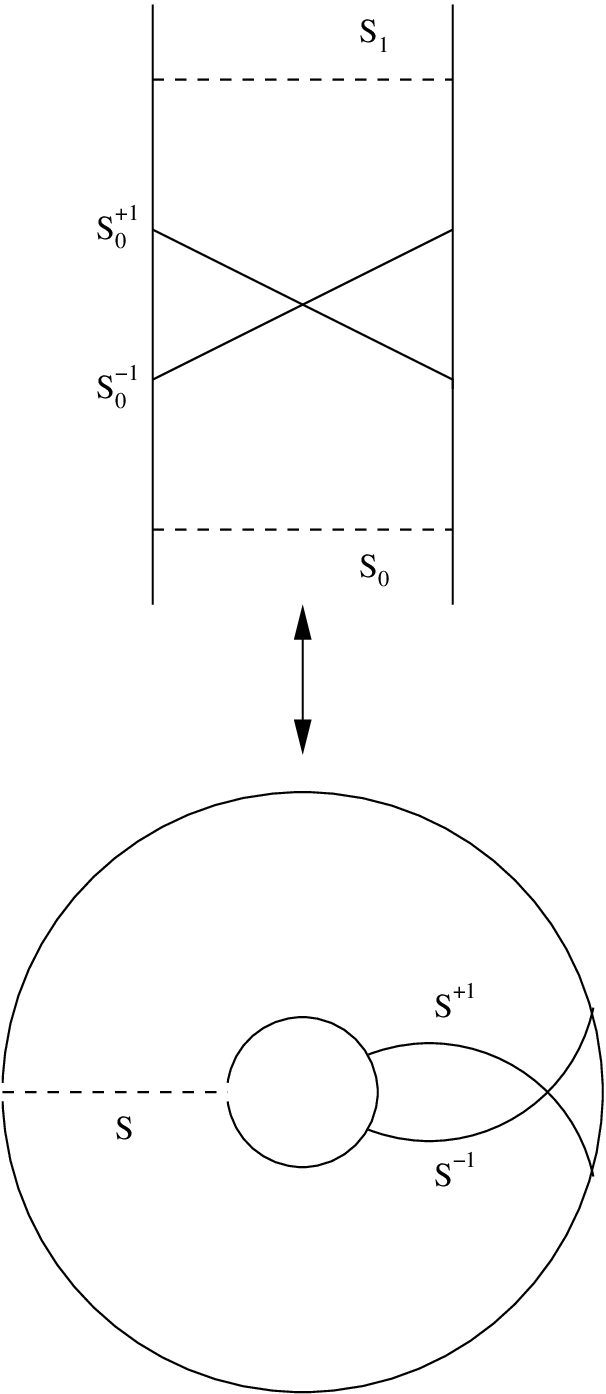} 
\caption{\sl $M(K)$ and $E(K)$}
\label{setup}
\end{figure}

\proof
We will construct $S^{up}, S^{down}$ explicitly by using the universal
abelian cover $M(K)$ of $E(K)$.  Let $\tau$ be a generator of the group of 
covering translations and let $S_0$ be a lift of $S$ to $M(K)$.  Set
\[S_1 = \tau(S_0).\]
Denote the component of \[M(K) \backslash (S_0 \cup S_1)\] that lies
between $S_0$ and $S_1$ by $C$ and note that $C$ is homeomorphic to
$E(K) \backslash S$ via the restriction of the covering map $M(K)\to
E(K)$.  See Figure \ref{setup}.  In particular, there are lifts
$S_0^{+1}, S_0^{-1}$ of $S^{+1}, S^{-1},$ respectively, in $C$.

Denote the two components of
$M(K) \backslash S_0^{-1}$ by $M^{-1}_-$ and $M^{-1}_+$, with $M^{-1}_+$
the component above $S_0^{-1}$, {\it i.e.}, the component containing $S_1$.
Similarly, denote the two components of 
$M(K) \backslash S^{+1}_0$ by $M^{+1}_-$ and $M^{+1}_+$, with $M^{+1}_+$
the component above $S^{+1}_0$, {\it i.e.}, the component containing $S_1$.
Recall that the {\em frontier} of a subset $H$, denoted by $fr(H)$, is the
closure of $H$ (in the ambient space) minus the interior of $H$.

\begin{figure}
[htbp] 
   \centering
    \includegraphics[width=3in]{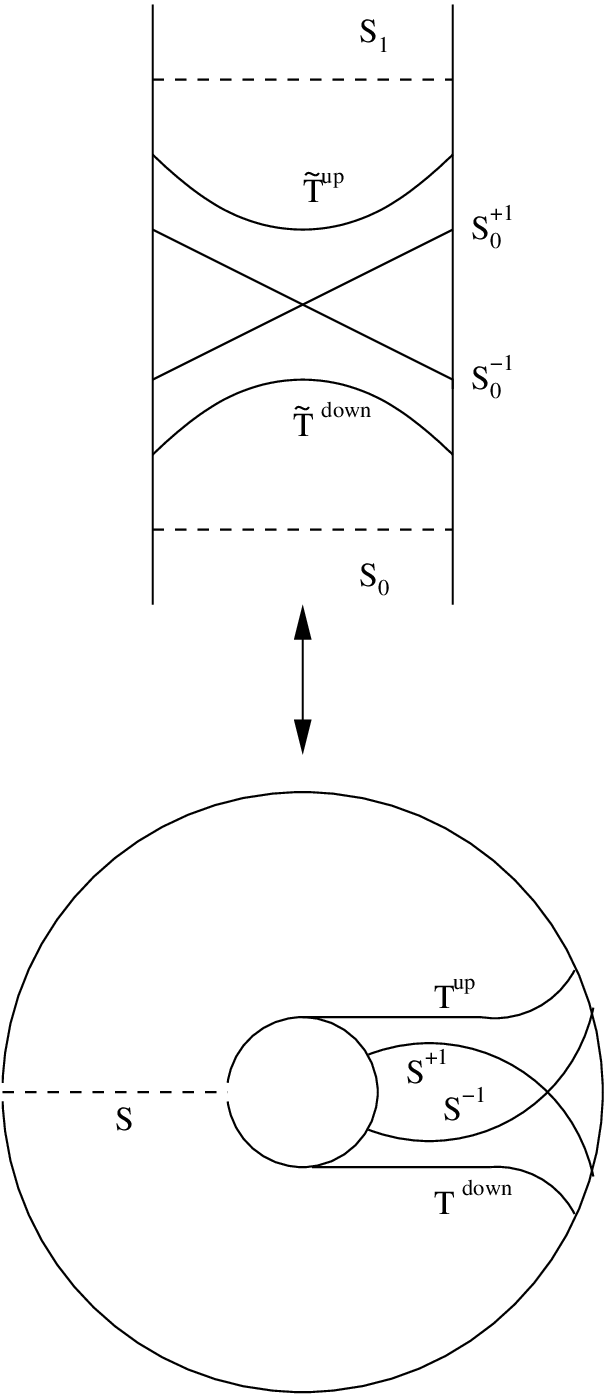} 
\caption{\sl $\tilde T^{down}$ and $\tilde T^{up}$}
\label{setup1}
\end{figure}

Set \[\tilde B^{down} = fr(M^{-1}_- \cap M^{+1}_-)\] and
\[\tilde B^{up} = fr(M^{-1}_+ \cap M^{+1}_+).\]
Further, let $\tilde T^{down}$ be a small pushoff of 
$\tilde B^{down}$ into $M^{-1}_- \cap M^{+1}_-$
and $\tilde T^{up}$ a small pushoff of $\tilde B^{up}$
into $M^{-1}_+ \cap M^{+1}_+$.  See Figure \ref{setup1}.
Let $B^{up}, B^{down}, T^{up}, T^{down}$ denote the (homeomorphic) projections of 
$\tilde B^{up}, \tilde B^{down}, \tilde T^{up}, \tilde T^{down}$ to the manifold $E(K)$. 
(Note that these surfaces could be disconnected.) 

\vspace{1 mm}
\noindent
Claim 1: $T^{up}, T^{down}$ each contain one Seifert surface.
\vspace{1 mm}

Let $\gamma$ be an oriented simple closed curve on $\partial E(K)$ that generates the 
homology of $E(K)$ and let $\tilde \gamma$ be its lift to $M(K)$.  
The algebraic intersection number of $\gamma$ with $\partial S, 
\partial S^{+1}, \partial S^{-1}$,
respectively, is $1$.  Hence the algebraic intersection number of $\tilde \gamma$
with $\partial S_0^{+1},$ $\partial S_0^{-1}$, respectively, is also $1$, because of the 
homeomorphism between $C$ and $E(K) \backslash S$.   
Note that one component of $\partial S_0^{- 1}\cup \partial S_0^{+1}$ lies in $\tilde B^{up}$ and the 
other lies in $\tilde B^{down}$.  It follows that the algebraic intersection number of $\tilde \gamma$ with 
$\partial \tilde T^{up}, \partial \tilde T^{down}$, respectively, is $1$.  
This in turn means that the algebraic
intersection number of $\gamma$ with $T^{up}, T^{down}$, respectively, is
$1$, because of the homeomorphism between $C$ and 
$E(K) \backslash S$.  It follows that $T^{up},$ $T^{down}$
each contain one Seifert
surface.  This proves Claim 1.

\vspace{1 mm}

Denote the Seifert surfaces in $T^{down}, T^{up}$ by $S^{down}, S^{up}$, respectively.
Then 
\[S \cap S^{down} = \emptyset\]\[S \cap S^{up} = \emptyset\]
 \[(S^{+1} \cup S^{-1}) \cap S^{down} = \emptyset\]
\[(S^{+1} \cup S^{-1}) \cap S^{up} = \emptyset\]
 \[S^{down} \cap S^{up} = \emptyset\]
 
 \vspace{1 mm}
 \noindent
 Claim 2: $S^{up}, S^{down}$ are minimal genus Seifert surfaces of $K$.
 \vspace{1 mm}

Note that
\[\chi(T^{up}) + \chi(T^{down}) = \chi(\tilde T^{up}) + \chi(\tilde T^{down})\]
since $T^{up} \cup T^{down}$ is homeomorphic to $\tilde T^{up} \cup \tilde T^{down}$.
Also, \[\chi(\tilde T^{up}) + \chi(\tilde T^{down}) = \chi(S_0^{+1}) + \chi(S_0^{-1}) ,\]
since $\tilde T^{up}$ (resp. $\tilde T^{down}$) is isotopic to $\tilde B^{up}$ 
(resp. $\tilde B^{down}$),  
and \[\tilde B^{up} \cup \tilde B^{down} = S_0^{+1} \cup S_0^{-1}.\]
Finally, \[\chi(S_0^{+1}) + \chi(S_0^{-1}) = \chi(S^{-1}) + \chi(S^{+1})\]
since $S_0^{-1} \cup S_0^{+1}$ is homeomorphic to $S^{-1} \cup S^{+1}$.
Therefore 
\[\chi(T^{up}) + \chi(T^{down}) = \chi(S^{-1}) + \chi(S^{+1}).\]
There are no disks in $S^{-1} \Delta S^{+1}$,
hence there are no disks in 
$S_0^{-1} \Delta S_0^{+1}$.  Thus there
are no disks or 2-spheres in $\tilde T^{up}, \tilde T^{down}$.  Therefore there are no
disks or 2-spheres in $T^{up}, T^{down}$.  
So
\[\chi(S^{up}) + \chi(S^{down}) \geq \chi(T^{up}) + \chi(T^{down})\]
since 
\[S^{up} \cup S^{down} \subset T^{up} \cup T^{down}.\]

Set $g = genus(K)$ (the minimal genus of a Seifert surface of $K$).  Then
\[\chi(S^{-1}) + \chi(S^{+1}) = 1 - 2g + 1 - 2g = 2 - 4g,\]
since $S^{+1}$ and $S^{-1}$ are minimal genus Seifert surfaces.  Note also that
\[\chi(S^{down}) \leq 1 - 2g, \chi(S^{up}) \leq 1 - 2g.\]  By this observation and the above computation,
\[2 - 4g \geq \chi(S^{up}) + \chi(S^{down}) \geq  2 - 4g.\]
It follows that \[\chi(S^{up}) = \chi(S^{down}) = 1 - 2g.\]
This proves Claim 2. 

\vspace{1 mm}

This proves part (a) of the Lemma.  Part (b) follows from the construction.
We now prove part (c) of the Lemma. 

Suppose that $E(K)$  is endowed with a 
relative Riemannian metric. Equip $M(K)$ with the pull-back of this metric. 
 Since $C$ is isometric to $E(K) \backslash S$ via the 
 restriction of the covering map $M(K)\to E(K)$, we have  
\[A(S^{-1}) + A(S^{+1}) = A(S_0^{-1}) + A(S_0^{+1})\] 
and
\[A(\tilde B^{down}) + A(\tilde B^{up}) = A(B^{down}) + A(B^{up}).\]
By the construction of $\tilde B^{down}, \tilde B^{up}$,  
\[A(S_0^{-1}) + A(S^{+1}) = A(\tilde B^{down}) + A(\tilde B^{up}).\]

Note that the surfaces $B^{down}, B^{up}$  are not smooth, 
while the relative least area surfaces in their respective isotopy classes are necessarily smooth. 
Therefore, 
$$
A(B^{down})> A([S^{down}]), \quad A(B^{up})> A([S^{up}]). 
$$
It follows that 
\[
A(S^{-1}) + A(S^{+1}) = A(\tilde B^{down})+ A(\tilde B^{up}) >  A([S^{down}]) + A([S^{up}]). 
\]
\qed

Lemma \ref{disj} below reinterprets Lemma \ref{disj1} above in the context 
of the Kakimizu complex.  

\begin{lem}  \label{disj}
Let $K$ be a knot in ${\mathbb S}^3$ and
suppose that $v, v^{+1}, v^{-1}$
are vertices in the Kakimizu complex of $K$ such that

\begin{enumerate}
\item \[d_K(v, v^{+1}) = 1, d_K(v, v^{-1}) = 1\]
\item \[d_K(v^{-1}, v^{+1}) = 2\] 
\item The complexity of $v$ is no smaller than the complexity of $v^{+1}, 
v^{-1},$ respectively.
\end{enumerate}

Then there are vertices $v^{down}, v^{up}$
in the Kakimizu complex of $K$ such that (a)
\[d_K(v, v^{down}) \leq 1, d_K(v, v^{up}) \leq 1, 
d_K(v^{+1}, v^{up}) \leq 1, d_K(v^{-1}, v^{up}) \leq 1, \]
\[d_K(v^{+1}, v^{down}) \leq 1, d_K(v^{-1}, v^{down}) \leq 1, 
d_K(v^{down}, v^{up}) \leq 1.\]
Moreover,  (b) if there is a vertex $w$ such that 
\[d_K(w, v^{+1}) = 1, d_K(w, v^{-1}) = 1\]
then \[d_K(w, v^{down}) \leq 1, d_K(w, v^{up}) \leq 1.\]
Furthermore, (c) the complexity of either $v^{up}$ or of $v^{down}$ is strictly 
less than that of $v$.
\end{lem}

\proof Endow $E(K)$ with a relative Riemannian metric and let ${\cal J}$ 
be a smooth foliation of $\partial E(K)$ by preferred longitudes. 
Let $S, S^{+1}, S^{-1}$ be relative least area representatives of 
$v, v^{+1}, v^{-1}$ which exist
by Theorem \ref{laexist}. 
We will say that $S^{+1}, S^{-1}$ are in {\em general position} 
if they intersect transversely and their boundaries in $\partial E(K)$ are disjoint.

\vspace{1 mm}
\noindent
Case 1: $S^{-1}, S^{+1}$ are in general position. 
\vspace{1 mm}

Since \[d_K(v, v^{+1}) = 1, d_K(v, v^{-1}) = 1,\]
\[d_K(v^{-1}, v^{+1}) = 2,\] 
we have \[S^{-1} \cap S^{+1} \neq \emptyset\] 
and Theorem \ref{ladisj} gives us 
 \[S \cap S^{+1} = \emptyset, S \cap S^{-1} = \emptyset.\] 

By \cite[Lemma 1.2]{FHS} there are no disks in $S^{-1} \Delta S^{+1}$. (Such a disk would yield a 
``product region'' in the sense of \cite{FHS}.) Thus all hypotheses of Lemma \ref{disj1} are satisfied.  
Let $v^{up}, v^{down}$ be the vertices
 in the Kakimizu of $K$
corresponding to $S^{up}, S^{down}$.  Then parts (a) and (b), respectively, follow from 
parts (a) and (b), respectively, of Lemma \ref{disj1}.

Furthermore, the statement about complexities follows because
say, \[A([S^{down}]) \geq A([S^{up}]),\]
and thus \[A(S) \geq max\{A(S^{-1}), A(S^{+1})\} \geq \frac{1}{2} (A(S^{-1}) + A(S^{+1}) )> \] 
\[\frac{1}{2}(A([S^{up}]) + A([S^{down}]))
\geq A([S^{up}]).\] Hence \[A(S) > A([S^{up}]).\]  This proves part (c).

\vspace{1 mm}
\noindent
Case 2: $S^{-1}, S^{+1}$ are not in general position. 
\vspace{1 mm}

In this case we apply the {\em Meeks-Yau trick} as described in \cite[Proof of Lemma 1.3]{FHS}: 
Let $x\in S^{-1}\cap S^{+1}$ be an interior point 
of  $E(K)$ where the surfaces intersect transversely. Let 
$D_{-1}\subset S^{-1}, D_{+1}\subset S^{+1}$ be small disks about $x$, 
both contained in the interior of $E(K)$, so that the intersection 
$\alpha=D_1\cap D_{-1}$ is a (smooth) arc of transverse intersection between these disks. 
Then there are portions of the surfaces $B^{up}, B^{down}$  near $x$ 
that are obtained by cutting $D_{+1}$ and $D_{-1}$ along $\alpha$ and pasting them together. 
The results are two piecewise-smooth disks $D^{up}\subset B^{up}$ and 
$D^{down}\subset B^{down}$. Neither disk is smooth along 
$\alpha$. Therefore, by ``rounding off'' these disks along $\alpha$ 
and keeping their boundaries fixed, we obtain two disks whose total area is less than 
$$
Area(D_1) + Area(D_2) -\epsilon, \quad \epsilon >0. 
$$ 
Next, take a surface  $S^{-1}(t)$ which is sufficiently close to $S^{-1}$ in the $C^1$-topology, so that:

1.  $\partial S^{-1}(t)\subset \partial E(K)$ lies in ${\cal J}$ and is disjoint from $\partial S^{+1}$.

2. $S^{-1}(t)$ intersects $S^{+1}$ transversely. 

3. $D_{-1}\subset S^{-1}(t)$. 

4. $Area(S^{-1}(t))< Area(S^{-1})+ \epsilon$. 
 
Now apply the argument from Case 1 to the surfaces $S^{-1}(t)$ and $S^{+1}$. As explained in \cite[Proof of Lemma 1.3]{FHS}, there are no ``product regions'' between the new surfaces. In particular, since  $S^{+1}$ and $S^{-1}$ are incompressible and $E(K)$ is irreducible, it follows that the symmetric difference $S^{-1}(t) \Delta S^{+1}$ contains no disks. Construct surfaces $B^{up}(t)$ and $B^{down}(t)$ and 
$S^{up}(t), S^{down}(t)$ in the same way as before. Since by  ``rounding-off'' $B^{up}(t)$ and $B^{down}(t)$  we loose more total area than we have gained by replacing $S^{-1}$ with $S^{-1}(t)$, we conclude that 
$$
A([S^{up}(t)])+ A([S^{down}(t)])< A(S^{-1})+ A(S^{+1}).  
$$  
The remainder of the argument follows as before. \qed

\section{The Kakimizu complex is flag}

We will use the following Theorem to prove the main theorem 
(Theorem \ref{kaksc}), but it is interesting in its own right.
Recall that a simplicial complex is {\it flag} if it contains no
empty simplices, {\it i.e.,} if it contains an $n$-simplex whenever it
contains the $(n-1)$-skeleton of the simplex.

\begin{thm} \label{flag}
The Kakimizu complex of a knot is flag.
\end{thm}

\proof
Let $K$ be a knot.  Endow $E(K)$ with a relative Riemannian metric.  
If the Kakimizu complex of $K$
contains the 1-skeleton of the simplex $\sigma$, then, by definition,
there are disjoint minimal genus Seifert surfaces representing any pair
of vertices in $\sigma$.  Hence, if we choose least area representatives 
for the vertices, it follows from Theorem \ref{ladisj} that these
representatives are simultaneously disjoint.  Thus
$\sigma$ belongs to the Kakimizu complex.
\qed

\section{The Kakimizu complex is simply connected}

We here prove that the Kakimizu complex is simply connected.  
Recall that paths and loops in the Kakimizu 
complex traverse only
vertices and edges (not higher dimensional simplices).  Recall
that we specify paths by the vertices they traverse, {\it e.g.},
$v^1, \dots, v^n$.    In the case of loops, we abuse notation slightly and 
write $0, 1, \dots, n+1$ when we really mean $0\; mod \; n, \dots, n+1\; mod \; n$.

\begin{thm} \label{kaksc}
  Let $K$ be a knot in ${\mathbb S}^3$.  The Kakimizu complex of $K$
  is simply connected.
\end{thm}

\proof Let $v^1, \dots, v^n$ be the vertices in a loop in the
Kakimizu complex of $K$.  

\vspace{1 mm}
\noindent
Claim 1: If $d_K(v^{i-1}, v^{i+1}) = 1$ for some $i$, then $v^1, \dots, v^n$
is homotopic to a shorter loop.  

If $d_K(v^{i-1}, v^{i+1}) = 1$ and $n = 3$, then the
loop spans a 2-simplex in the Kakimizu complex by Theorem \ref{flag} 
and is hence homotopic to a single vertex.
If $d_K(v^{i-1}, v^{i+1}) = 1$ and $n > 3$, then
$v^{i-1}, v^i, v^{i+1}$ still spans a 2-simplex.  Hence the loop
$v^1, \dots, v^n$ is homotopic to the loop obtained by
replacing $v^{i-1}, v^i, v^{i+1}$ with $v^{i-1},
v^{i+1}$.  This proves Claim 1.

\vspace{1 mm}

We will henceforth assume that
\[d_K(v^{i-1}, v^{i+1}) = 2 \;\; \forall i .\]

\vspace{1 mm}
\noindent
Claim 2: If $c(v^i) \geq max\{c(v^{i-1}), c(v^{i+1})\}$, 
then the complexity of the loop $v^1, \dots, v^n$ is not minimal
among loops homotopic to $v^1, \dots, v^n$.
 
Under the assumption $c(v^i) \geq max\{c(v^{i-1}), c(v^{i+1})\}$, Lemma
\ref{disj} furnishes two vertices $v^{up},$ $v^{down}$ in the Kakimizu complex such that
\[d_K(v^{i+1}, v^{up}) \leq 1, d_K(v^{i-1}, v^{up}) \leq 1,
d_K(v^{i+1}, v^{down}) \leq 1, d_K(v^{i-1}, v^{down}) \leq 1,\]
\[d_K(v^i, v^{down}) \leq 1, d_K(v^i, v^{up}) \leq 1\]
and such that the complexity of say, $v^{up},$ 
is strictly less than that of $v^i$.

\begin{figure}[htbp] 
   \centering
    \includegraphics[width=3in]{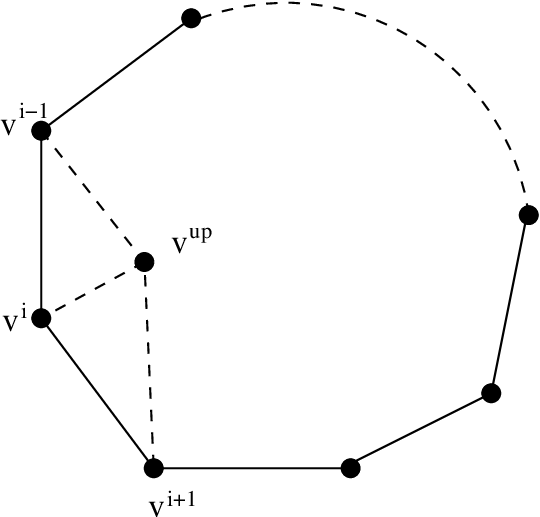} 
\caption{A homotopy}
\label{loop}
\end{figure}

In particular, the two loops 
\[v^{i-1}, v^i, v^{up}\] \[v^i, v^{i+1}, v^{up}\] bound 
2-simplices in the Kakimizu complex of $K$, so 
\[v^1, \dots, v^n\] is homotopic to 
\[v^1, \dots, v^{i-1}, v^{up}, v^{i+1}, \dots, v^n\]
See Figure \ref{loop}.  

Since $c(v^{up}) < c(v^i)$, the complexity of this path is strictly smaller than that of
the path $v^1, \dots, v^n$.  This proves Claim 2.

\vspace{1 mm}

Now suppose that there is a homotopically nontrivial loop in the Kakimizu complex of $K$.
By Theorem \ref{min}, there is a homotopically nontrivial loop 
$v^1, \dots, v^n$ of smallest complexity.  Then Claim 2 tells us that
\[c(v^i) < max\{c(v^{i-1}), c(v^{i+1})\} \; \forall i.\]
But this is impossible.  
Thus the Kakimizu
complex of $K$ is simply connected.  \qed

\section{A note on contractibility} \label{contractibility}

We recall some of the standard terminology for simplicial
complexes: The {\em link} of a vertex $v$ in a simplicial
complex $X$, denoted by $X_{v}$, is the union of all
simplices disjoint from $v$ that together with $v$ span
a simplex in $X$. The {\em star} of a vertex $v$ in a
simplicial complex $X$ is the union
of all simplices in $X$ that contain $v$.  
A 2-dimensional simplicial complex $X$ is said to be 
{\em locally $k$-large} if for every vertex $v\in X$, every 
homotopically nontrivial loop in $X_v$ has 
length at least $k$. 

The following theorem is classical.  It follows from Propositions
II.4.1 (Cartan-Hadamard Theorem) and II.5.25 in \cite{BH}.  See
Section \ref{speculation} for a more general version.

\begin{thm} \label{classical}
  The universal cover of a 2-dimensional connected locally
  6-large simplicial complex is contractible.
\end{thm}

\begin{lem} \label{links1}
  A homotopically nontrivial loop in the link of a vertex in the
  Kakimizu complex must have length at least $5$.
\end{lem}

\proof Loops of length 3 in the Kakimizu complex bound 2-simplices,
by Theorem \ref{flag}.
Hence it suffices to show that there are no homotopically nontrivial
loops of length 4 in the link of a vertex.  Suppose that 
$v^1, v^2, v^3, v^4$ is a loop in
the link, $X_v$, of the vertex $v$.  We argue as in the proof of
Theorem \ref{kaksc}.  The same reasoning as in Claim 1 of the proof
of Theorem \ref{kaksc} tells us that we need only consider the
case in which \[d_K(v^{i-1}, v^{i+1}) = 2 \;\;\forall i.\] 
The reasoning in Claim 2 of the proof of Theorem \ref{kaksc},
tells us that there is a loop homotopic to $v^1, \dots, v^4$ of lesser
complexity.  By part (b) of Lemma \ref{disj}, the vertex $v^{up}$ that replaces
one of $v^1, \dots, v^4$ lies in the star of $v$.  It either lies in the link of $v$ or it  
is equal to $v$.  In the latter case we must examine how this comes about.  Following
the reasoning in Claim 2 of the proof of Theorem \ref{kaksc}, we
apply Lemma \ref{disj} to the vertices $v^i, v^{i+1}, v^{i-1}$ in place of 
the vertices $v, v^{+1}, v^{-1}$ to obtain vertices $v^{down}, v^{up}$ such
that 

\[d_K(v^{i+1}, v^{down}) \leq 1, d_K(v^{i-1}, v^{down}) \leq 1, 
d_K(v^{i+1}, v^{up}) \leq 1, d_K(v^{i-1}, v^{up}) \leq 1, \]
\[d_K(v^i, v^{down}) \leq 1, d_K(v^i, v^{up}) \leq 1, d_K(v, v^{down}) \leq 1, 
d_K(v, v^{up}) \leq 1\]
and (by part (b) of Lemma \ref{disj}) also that
\[d_K(v^{i+2}, v^{down}) \leq 1, d_K(v^{i+2}, v^{up}) \leq 1.\] 

Two cases need to be considered:

\vspace{1 mm}
\noindent
Case 1: $v^{up} = v$ and $v^{down} \neq v$.

\vspace{1 mm}
  
\begin{figure}[htbp] 
   \centering
    \includegraphics[width=2in]{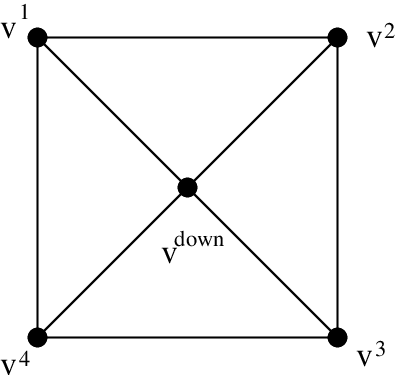} 
\caption{A homotopically trivial loop}
\label{5largea}
\end{figure}

If $v^{down}$ is equal to one of the $v^j$, then the loop $v^1, v^2, v^3, v^4$ breaks into
two loops of length 3.  The latter bound 2-simplices by Theorem \ref{flag}.  Hence we will
assume that $v^{down}$ is not equal to one of the $v^j$.  Now 
$v^{i-1}, v^i, v^{down}$ and $v^i, v^{down}, v^{i+1}$ are loops of length 3 
and hence bound 2-simplices by Theorem \ref{flag}.  It follows that our loop 
\[v^{i-1}, v^i, v^{i+1}, v^{i+2}\] is homotopic 
to \[v^{i-1}, v^{down}, v^{i+1}, v^{i+2}\]
in $X_v$.  In addition, both $v^{down}, v^{i+1}, v^{i+2}$ 
and  $v^{i+2}, v^{i-1}, v^{down}$ bound 2-simplices by Lemma \ref{flag}.  
Thus $v^{i-1}, v^{down}, v^{i+1}, v^{i+2}$ and hence $v^{i-1}, v^i, v^{i+1}, v^{i+2}$,
{\it i.e.,} $v^1, \dots, v^4$, is homotopically trivial in $X_v$.  See Figure \ref{5largea}.

\vspace{1 mm}
\noindent
Case 2:  $v^{up} = v^{down} = v.$ 

\vspace{1 mm}

\begin{figure} [htbp] 
   \centering
    \includegraphics[width=2in]{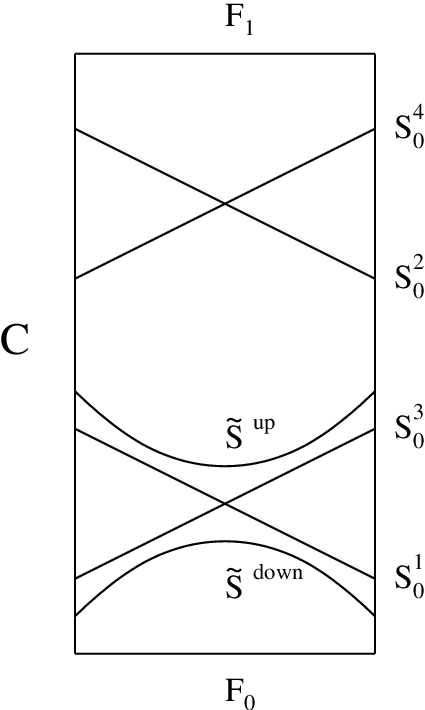} 
\caption{The product region}
\label{5largeb}
\end{figure}

In this case we must argue differently.  Let $F$ be a representative of $v$,
$F_0$ a lift of $F$ to $M(K),$ $F_1 = \tau(F_0)$ and $C$ the component of 
\[M(K) \backslash (F_0 \cup F_1)\] that lies between $F_0$ and $F_1$.
Recall the notation $S_0^j,$ 
$\tilde T^{up},$ $\tilde T^{down},$ $S^{up},$ 
$S^{down}$ from the proof of Lemma \ref{disj1}
and let $\tilde S^{up}, \tilde S^{down}$ be lifts of $S^{up}, S^{down}$ to $C$, respectively.
See Figure \ref{5largeb}.  We will assume (only) that $F,$ $F_0,$ $F_1$ and each $S_0^j$ are 
least area surfaces in their isotopy classes.  

Denote the components of $M(K) \backslash S_0^j$ by
$M_-^j$ and $M_+^j$, with $M_+^j$ the component above $S_0^j,$ {\it i.e.,} the component
containing $F_1$.  (Note: we say that a surface in $C$ lies {\it above} $S_0^j$ if
it lies in $M_+^j$ and {\it below} $S_0^j$ if it lies in $M_-^j$.  It is easy to check,
but important to realize, that if $S_0^j$ lies above $S_0^k$, then $S_0^k$ lies 
below $S_0^j$.)
It follows from the construction of $S^{up},$ 
$S^{down}$ in Lemma \ref{disj1} that $\tilde S^{up},$ $\tilde S^{down}$ are contained 
in $\tilde T^{up},$ $\tilde T^{down},$ respectively.   Moreover, the push-off yielding 
$\tilde T^{up}$ forces it to lie above $S_0^{i-1}$ and $S_0^{i+1},$ 
{\em i.e.}, in the interior of $M_+^{i-1}$ and $M_+^{i+1}$.  
The push-off yielding $\tilde T^{down}$
forces it to lie below $S_0^{i-1}$ and $S_0^{i+1},$ {\em i.e.}, 
in the interior of
$M_-^{i-1}$ and $M_-^{i+1}$.  

Here $S_0^i,$ $S_0^{i+2}$ are disjoint
from $S_0^{i-1} \cup S_0^{i+1}$ by Theorem \ref{ladisj} and hence, after isotopy if necessary, also from
$\tilde S^{up} \cup \tilde S^{down}$.  
We wish to show that $S_0^i, S_0^{i+2}$ either both lie below $\tilde S^{down}$ 
or both lie above $\tilde S^{up}$.   Indeed, since $S_0^i \cap S_0^{i+2} \neq \emptyset,$
we need only exclude the possibility that $S_0^i \cup S_0^{i+2}$ lie above
$\tilde S^{down}$ but below $\tilde S^{up}$.  Consider how $S_0^{i-1} \cup S_0^{i+1}$
lie with respect to $S_0^i$.  Since $S_0^{i-1} \cap S_0^{i+1} \neq \emptyset,$ they
either both lie in $M_-^i$ or both lie in $M_+^i$.  If \[S_0^{i-1} \cup S_0^{i+1} \subset M_-^i,\] 
then $S_0^i$ lies above $S_0^{i-1} \cup S_0^{i+1}$ and hence above $\tilde S^{up}$.  
If \[S_0^{i-1} \cup S_0^{i+1} \subset M_+^i,\] 
then $S_0^i$ lies below $S_0^{i-1} \cup S_0^{i+1}$ and hence below $\tilde S^{down}$.
A similar argument applies to $S_0^{i+2}$.  Thus $S_0^i, S_0^{i+2}$ either both lie 
below $\tilde S^{down}$ or both lie above $\tilde S^{up}$.

Since $v^{up} = v^{down} = v$, $\tilde S^{up}$ is isotopic
to $F_0$ or $F_1$ as is $\tilde S^{down}$.  
It follows that either $S^{i-1}$ and $S^{i+1}$ lie in a
product region and are hence isotopic, or that $S^{i}$ and $S^{i+2}$
lie in a product region and are hence isotopic.  See Figure \ref{5largeb}.
But this is
impossible.  Hence there are no homotopically nontrivial loops of
length less than or equal to $4$ in the link of $v$.  \qed

\begin{lem} \label{links2}
  A homotopically nontrivial loop in the link of a vertex in the 
  Kakimizu complex must have length at least $6$.
\end{lem}

\proof
By Lemma \ref{links1}, it suffices to show that there are no
homotopically nontrivial loops of length 5.
So suppose that $v^1, \dots, v^5$ is a loop of length $5$ in
$X_v$.  By the reasoning in Claim 1 in the proof of Theorem \ref{kaksc}, 
we may assume that $d_K(v^{i-1}, v^{i+1}) = 2$ $\forall i$. 
(For otherwise the loop is homotopic, in $X_v$, to a loop of length
$4$, and hence homotopically trivial by Lemma \ref{links1}.)
Endow $E(K)$ with a relative Riemannian metric and 
let $S^{1}, \dots, S^5$ be relative least area representatives
of $v^1, \dots, v^5$.  By Theorem \ref{ladisj},
\[S^i \cap (S^{i+1} \cup S^{i-1}) = \emptyset \;\;\ \forall i\]
and \[S^{i-1} \cap S^{i+1} \neq \emptyset \;\; \forall i.\]
 
\begin{figure}[htbp] 
   \centering
    \includegraphics[width=2in]{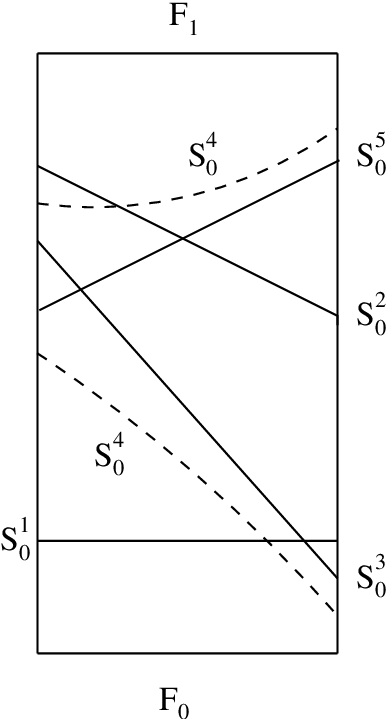} 
\caption{Five separating surfaces}
\label{6large}
\end{figure}

Recall the notation $M(K), F_0, F_1, C$ from the proof of Lemma \ref{links1}
and let $S_0^1, \dots, S_0^5$ be lifts of $S^{1}, \dots, S^5$
to $C$.
Then
\[S^i_0 \cap (S^{i+1}_0 \cup S^{i-1}_0) = \emptyset \; \; \forall i\] and
\[S^i_0 \cap S^j_0 \neq \emptyset \; {\rm for} \; j \neq i, i \pm 1.\]
Each $S^i_0$ is separating in $C$.  In
particular, $S_0^2 \cap S_0^5 \neq \emptyset,$ so these two surfaces
must either both lie above or both lie below $S_0^1$.  Assume the
former, as the other case will then follow by a symmetric argument.  Also,
$S_0^3$ must lie below $S_0^2$, in order to have nonempty
intersection with $S_0^1$. Furthermore, $S_0^4$ must lie above
$S_0^3$ in order to have nontrivial intersection with $S_0^2$, but
below $S_0^5$, in order to have nontrivial intersection with
$S_0^1$.  See Figure \ref{6large}.

Since $S_0^3 \cap S_0^5 \neq \emptyset$, this is
impossible. \qed

\begin{thm} \label{contractible}
  If the Kakimizu complex of a knot $K$ is at most 2-dimensional, then
  it is contractible.
\end{thm}

\proof Lemmas \ref{links1} and \ref{links2} establish the fact that links of 
vertices in the Kakimizu complex 
contain no nontrivial loops of lengths $\le 5$. Thus the hypotheses of 
Theorem \ref{classical} are satisfied,
whence the universal cover of the Kakimizu complex is contractible.
Hence by Theorem \ref{kaksc} the Kakimizu complex is contractible.
\qed

\section{Further remarks} \label{speculation}

To improve readability of this paper, the results have been stated 
and proved for knots in ${\mathbb S}^3$.  However, they are equally 
valid for knots in homology 3-spheres.  
Furthermore, the Kakimizu complex is defined in terms of isotopy classes of
minimal genus Seifert surfaces.  A related complex is obtained by
considering all isotopy classes of incompressible Seifert
surfaces.  The theorems and arguments used here are stated in
terms of minimal genus Seifert surfaces.  It seem likely that this assumption is
unnecessary and that the analogous theorems hold for the more
general complex as well.

The goal of proving contractibility of the Kakimizu complex may be
out of reach, but there are natural questions to ask, now that
simple connectedness is established. One such question is whether
or not the Kakimizu complex is 2-connected. Indeed, this question
is being pursued by Sakuma and Shackleton (\cite{Sh}) who believe
that it is 2-connected and that this can be established via the
techniques used here.

Theorem \ref{classical} has a generalization, due to Januszkiewicz and \'Swi{\c a}tkowski (\cite{JS}), 
to higher-dimensional simplicial complexes, although their notion of local $k$-largeness is more subtle. 
It is unclear if Lemmas \ref{links1} and \ref{links2} can be 
modified to fit their definition.

\section{Appendix: A compactness theorem for stable minimal surfaces \\ by Michael Kapovich}

Let $M$ be a ${\mathbb P}^2$-irreducible compact Riemannian 3-manifold with smooth strictly convex boundary and 
${\cal J}$ a compact family of smooth curves on $\D M$. In the setting here we are mostly interested in the case where 
$\D M$ is a single torus and ${\mathcal J}$ is a smooth foliation of $\D M$ by closed curves. 
Let $S$ be a compact connected surface, possibly with boundary. 

Given a smooth proper embedding $f: (S, \D S)\to (M, \D M)$ we let $[f]$ denote its (proper) isotopy class (here we are not fixing the boundary value $f|\D S$). From now on, we will assume that the isotopy classes $[f]$ are such that for each boundary component $\D_i M$ of $M$, $f^{-1}(\D_i M)$ is a single component of $\D S$ and that the surface $f(S)$ is incompressible. It follows that for each  $f$ as above, there exists a least area (parameterized) surface in the set
$$
\{g\in [f]: g| \D S=f|\D S\}, 
$$
see \cite[Theorem 6.12]{Hass-Scott(1988)}. Such a surface is necessarily a stable minimal surface.  

Define the ``moduli space'' $\M([f])$ of stable minimal surfaces in the given proper isotopy class $[f]$ subject to the condition that 
$f|\D S$ is a parameterized multi-curve  in ${\cal J}$. 
Here we identify parameterized surfaces which differ by a reparameterization of $S$. We let $\M=\M(S)=\cup_{[f]} \M([f])$ be the space of all stable embedded minimal surfaces of the given topological type.  We give $\M$ the $C^1$-topology. 
Given a number $a$ we let $\M_a$ denote the subset of $\M$ consisting of surfaces of area $\le a$. 

\begin{prop}\label{comp}
If $S$ has nonempty boundary, then the space $\M_a$ is compact.  
\end{prop}

When $\D S$ is empty, a similar compactness result holds by a theorem of Nakauchi.  In \cite{Nakauchi}, Nakauchi uses Schoen's estimates \cite{Schoen(1983)} on the norm of the 2-nd fundamental form of stable minimal surfaces away from the boundary of $M$ to conclude that sequences of stable minimal surfaces admit convergent subsequences, except that the limiting surfaces in \cite{Nakauchi} may fail to be embedded but appear as 2-fold coverings of embedded surfaces. 

\proof
In the case of surfaces with boundary we modify Nakauchi's argument as follows: Given a sequence of minimal surfaces $f_i: S\to M$ whose boundary values are in ${\cal J}$ and whose area is $\le a$, a theorem of Anderson \cite[Theorem 3.1]{Anderson(1985b)} tells us that there exists a subsequence $f_{i_j}$ so that the sequence $f_{i_j}(S)$ converges to a minimal surface $\Si\subset M$ whose boundary is in ${\cal J}$. The surface $\Si$ need not be of the same topological type as $S$. However, the convergence of the surfaces is smooth away from a finite subset $x_1,...,x_m\in \Si$. Moreover, Anderson proves \cite[Paragraph 4 of the proof of Theorem 3.1]{Anderson(1985b)} that all the points $x_i$ belong to the interior of $M$. 
Thus we can apply Schoen's estimates  \cite{Schoen(1983)} to each point $x_i\in int(M)$ 
in the same manner as Nakauchi does, provided that the $f_i(S)$ are stable minimal surfaces. Schoen's estimates imply smooth convergence at the points $x_1,...,x_m$. Therefore, the maps $f_{i_j}: S\to M$ converge smoothly to a covering map $f: S\to \Si$. Suppose that $f$ is a nontrivial covering, then 
its restriction to $\D S$ is also nontrivial. However, compactness of ${\cal J}$ implies that the maps $f_{i_j}: \D S \to \D M$ converge to an embedding. This is a contradiction. Therefore, $f$ is 1-1. \qed 
 
\begin{rem}
If $M$ is a closed 3-manifold with a bumpy Riemannian metric, Colding and Minicozzi proved in  
\cite{Colding-Minicozzi(2000)} that the space of (not necessarily stable) minimal surfaces of uniformly bounded area is finite: Bumpiness of the metric is used to ensure that the limiting minimal surface $\Si$ has no nontrivial Jacobi fields. The same argument can be used in conjunction with the above proposition to ensure finiteness of $\M_a$, provided that the metric on $M$ is chosen to be bumpy on the interior of the manifold. 
\end{rem}

\begin{cor}\label{finite}
$\M_a$ splits as a disjoint union of finitely many open and closed subsets $\M_a([g])$ each of which consists of isotopic surfaces. 
\end{cor}
\proof Let $g$ be a stable minimal surface in  $\M_a$ and let $f_k: S\to M$ be a 
sequence of surfaces in $\M_a \backslash \M_a([g])$.  Suppose that the sequence $f_k: S\to M$ converges to
$f: S\to M$. Then $f$ is isotopic to each $f_k$ for $k$ sufficiently large. In particular, $f$ cannot be in $\M_a([g])$.
It follows that the subsets $\M_a([g])$ are open.  
Finiteness of the number of these sets follows from the compactness of $\M_a$. \qed 

\begin{cor} \label{exists}
Each $\M([f])$ contains an area-minimizer.  In particular, each class $\{g\in [f]: g|\partial S \in {\mathcal J}\}$ admits an area-minimizer. 
\end{cor}
\proof Area is a continuous function, therefore the corollary follows from the compactness of $\M_a([f])$. \qed 

\medskip
We now assume that $\partial M$ is a union of tori. 
We say that a Riemannian metric on $M$ is {\em relative} if $M$ 
has a strictly convex boundary. It is easy to see that $M$ always 
admits a relative Riemannian metric: Start with a metric on $T^2\times [0,1]$, 
where $T^2\times \{0\}$ strictly convex. Identify  $T^2\times [0,1]$ with neighborhoods of boundary tori 
$\partial_i M$ of $M$, where  $\partial_i M$ 
is identified with $T^2\times \{0\}$. Then extend the metric on a neighborhood of $\partial M$ 
arbitrarily to the rest of $M$. Fix a $C^1$-foliation ${\cal J}$ of $\partial M$ by closed curves. 

The area minimizers in the isotopy classes $\M([f])$ will be called 
{\em relative least area surfaces} (the word ``relative'' refers to the 
fact that we are assuming that $f(\partial S)\in {\cal J}$). 
We now compare our setup with that of \cite{FHS}. 

In our boundary conditions for relative least area surfaces we are using boundary curves in ${\cal J}$ rather than using a free boundary condition as is done in \cite{FHS}. However, since curves in $ {\cal J}$ are either disjoint or equal, this ensures that the cut-and-paste arguments used in \cite{FHS} preserve our boundary conditions. With this in mind,  the arguments in \cite{FHS} go through with our set-up. In particular, we obtain:

\begin{thm}\label{disjoint}
Let $f_i: S_i\to M, i=1,...,n,$ be incompressible surfaces which are pairwise non-isotopic and pairwise disjoint. 
Let $g_i: S_i\to M, i=1,...,n$ be relative area minimizers in the isotopy classes of $f_i, i=1,...,n$. Then $g_1(S_1),...,g_n(S_n)$ are also 
pairwise disjoint. 
\end{thm} 

We now assume that $M$ is diffeomorphic to the exterior of a knot in 
${\mathbb S}^3$, ${\cal J}$ is a foliation of $\partial M$ by preferred longitudes
and the surface $S$ has a single boundary component. 

Fix a number $n$ and let $\L \subset \M^n$ denote the subset of $n$-tuples 
represented by parameterized minimal genus Seifert surfaces $(f_1,...,f_n)$, so that:

a)  $d_K([f_i(S)], [f_{i+1}(S)])= 1$, $i$ is taken mod $n$. 
Here $d_K$ is the Kakimizu distance and the isotopy class $[f_k(S)]=v_k$ represents a  vertex in the Kakimizu complex. 

b) The loop in the 1-skeleton of the Kazimizu complex represented by the vertices $v_1,...,v_n$ is 
homotopically nontrivial in the Kazimizu complex. 

It follows immediately from Corollary \ref{finite} that $\L$ is closed in $\M^n$.   
Define the area functional $A: \M^n  \to \R_+$,  
$$
A(f_1,...,f_n)=\sum_{i=1}^n A(f_i(S)). 
$$

We then obtain the following:

\begin{cor} \label{mini}
The functional $A|\L$ attains a minimum. 
\end{cor}
\proof It is clear that $A$ is continuous and positive. By Proposition \ref{comp}, $A$ is proper. Since $\L$ is closed in 
$\M^n$, the restriction $A|\L$ is proper as well. Therefore, $A|\L$ attains its minimum. \qed

\vspace{2 mm}

\noindent
Department of Mathematics

\noindent
1 Shields Avenue

\noindent
University of
California, Davis

\noindent
Davis, CA 95616

\noindent
USA

\end{document}